\documentclass[letterpaper, 10 pt, conference]{ieeeconf}  

\IEEEoverridecommandlockouts                              

\usepackage{amsmath} 	
\usepackage{amssymb}  	
\usepackage{amsfonts}
\usepackage{graphicx}
\usepackage{subfigure}
\usepackage{epstopdf}
\usepackage{dblfloatfix}
\usepackage{tikz}
\usetikzlibrary{shapes,matrix,decorations.shapes}
\usepackage{array,xcolor}
\usepackage{algorithm}
\usepackage{algpseudocode}
\usepackage{multirow}
\usepackage{threeparttable}
\usepackage{booktabs}
\usepackage{balance}
\usepackage{cite}

\newtheorem{definition}{Definition}
\newtheorem{theorem}{Theorem}

\newtheorem{remark}{Remark}

\newtheorem{corollary}{Corollary}

\usepackage[colorlinks=true,      
linkcolor=black,      
citecolor=black,      
filecolor=black,      
urlcolor=blue ]{hyperref}

\title{\LARGE \bf
Improving Efficiency and Scalability of Sum of Squares Optimization: Recent Advances and Limitations\thanks{This is an invited tutorial paper for the 2017 IEEE International Conference on Decision and Control. Most details are omitted and can be found in the relevant references.}
}

\author{Amir Ali Ahmadi$^{1}$, Georgina Hall$^{1}$, Antonis Papachristodoulou$^{2}$, James Saunderson$^{3}$, and Yang Zheng$^{2}$
	\thanks{$^{1}$A. A. Ahmadi and G. Hall are with the department of Operations Research and Financial Engineering,
		Princeton University, Princeton, NJ 08540, USA.
		Emails: {\tt\small a\_a\_a@princeton.edu}; {\tt\small gh4@princeton.edu}}%
	\thanks{$^{2}$Y. Zheng and A. Papachristodoulou are with the Department of Engineering Science, University of Oxford, Parks Road, Oxford OX1 3PJ, U.K. Emails: {\tt \small yang.zheng@eng.ox.ac.uk}; {\tt \small antonis@eng.ox.ac.uk} }
	\thanks{$^{3}$J. Saunderson is with the Department of Electrical and Computer Systems Engineering,
		Monash University, VIC 3800, Australia. Email: {\tt \small james.saunderson@monash.edu}}%
}

\begin{document}

	\maketitle
	\thispagestyle{empty}
	\pagestyle{empty}

\begin{abstract}

It is well-known that any sum of squares (SOS) program can be cast as a semidefinite program (SDP) of a particular structure and that therein lies the computational bottleneck for SOS programs, as the SDPs generated by this procedure are large and costly to solve when the polynomials involved in the SOS programs have a large number of variables and degree. In this paper, we review SOS optimization techniques and present two new methods for improving their computational efficiency. The first method leverages the sparsity of the underlying SDP to obtain computational speed-ups. Further improvements can be obtained if the coefficients of the polynomials that describe the problem have a particular sparsity pattern, called \emph{chordal sparsity}. The second method bypasses semidefinite programming altogether and relies instead on solving a sequence of more tractable convex programs, namely linear and second order cone programs. This opens up the question as to how well one can approximate the cone of SOS polynomials by second order representable cones. In the last part of the paper, we present some recent negative results related to this question.
\end{abstract}

\section{Introduction and sum of squares review} \label{sec:intro}

Polynomial optimization is the problem of minimizing a (multivariate) polynomial function on a basic semialgebraic set; \emph{i.e.}, a subset of the Euclidean space defined by polynomial equations and inequalities.
This is an extremely broad class of optimization problems, with high-impact application areas throughout engineering and applied mathematics, which until not too long ago was believed to be hopelessly intractable to solve computationally. In recent years, however, a fundamental and exciting interplay between \emph{convex optimization} and \emph{algebraic geometry} has allowed for the solution or approximation of a large class of (nonconvex) polynomial optimization problems.

Amazingly, the success of this area stems from the ability to work around a single central question which is very simple to state: how can one test if a polynomial
$$
    p(x)\mathrel{\mathop:}=p(x_1,\ldots,x_n)
$$
is \emph{nonnegative}, \emph{i.e.}, satisfies $p(x)\geq 0$ for all $x\in\mathbb{R}^n$?

Unfortunately, answering this question is NP-hard already when $p(x)$ has degree 4. A powerful and more tractable sufficient condition for nonnegativity of $p(x)$, however, is for it to be a \emph{sum of squares} polynomial. A sum of squares (SOS) polynomial $p(x)$ is a polynomial that can be written as
{\color{black}
$$
    p(x)=\sum_{i=1}^r f_i^2(x)
$$
for some other polynomials $f_i(x), i = 1, \ldots, r$}. The question as to whether a nonnegative polynomial can always be written as a sum of squares has a celebrated history, dating back to Hilbert's 17th problem~\cite{Reznick} around the year 1900. What has caused a lot of recent excitement, however, is the discovery that the task of testing the SOS property and finding a sum of squares decomposition can be fully automated. This is a consequence of the following characterization of the set of SOS polynomials: a polynomial $p(x)$ of degree $2d$ is SOS if there exists a positive semidefinite matrix $Q$ (usually called the \emph{Gram matrix} of $p$) such that
\begin{align} \label{eq:charac.sos.sdp}
    p(x)=z(x)^TQz(x),
\end{align}
where $z(x)=[1,x_1,\ldots,x_n,\ldots,x_n^d]$ is the standard vector of monomials of degree $d$~\cite{sdprelax}. Hence, testing whether a polynomial is a sum of squares amounts to solving a \emph{semidefinite program} (SDP), a class of convex optimization problem for which numerical solution methods are available.

This simple but fundamental discovery forms the basis of a modern subfield of mathematical programming called \emph{``sum of squares optimization''}. An \emph{SOS program} is an optimization problem of the following form:

\begin{equation}\label{eq:SOS.programs}
\begin{aligned}
&\min_{p} &&C(p)\\
&\text{s.t. } &&A(p)=b\\
& &&p \text{ is SOS},
\end{aligned}
\end{equation}
where the decision variables are the coefficients of the polynomial $p$, the objective $C(p)$ is some linear function of the coefficients of $p$, and $A(p)$ are affine constraints in the coefficients of $p$. As a consequence of the aforementioned characterization of SOS polynomials, this program can be recast as the following semidefinite program:
\begin{equation}\label{eq:SOS.SDP}
\begin{aligned}
&\min_{p,Q} &&C(p)\\
&\text{s.t. } &&A(p)=b\\
& &&p(x)=z(x)^TQz(x),~\forall x\\
& &&Q\succeq 0.
\end{aligned}
\end{equation}

The most direct consequence of SOS optimization is for polynomial optimization: Under mild assumptions, the (global) minimum of a polynomial $p(x)$ on a basic algebraic set $K$ turns out to be equal to the largest scalar $\gamma$ such that $p(x)-\gamma$ is certified to be nonnegative on $K$ with a sum of squares proof~\cite{sdprelax},~\cite{Lglob2001}. The power of this statement stems from the fact that \emph{no convexity assumption} is placed on the polynomial optimization problem and yet the search for the sum of squares certificates is a convex (in fact semidefinite) program.

Aside from polynomial optimization, numerous other areas of computational mathematics have been impacted by sum of squares techniques. These include approximation algorithms for NP-hard combinatorial optimization problems~\cite{Stability_number_SOS},~\cite{barak2014sum}, equilibrium analysis of games~\cite{Pablo_poly_games}, robust and stochastic optimization~\cite{DanIancu_sos}, statistics and machine learning~\cite{barak2014dictionary},~\cite{convex_fitting}, software verification~\cite{MardavijRoozbehani2008},~\cite{roozbehani2005convex}, filter design~\cite{filter_design_trig_sos}, quantum computation~\cite{Pablo_Sep_Entang_States}, automated theorem proving~\cite{harrison2007verifying}, and fault diagnosis and verification of hypersonic aircraft~\cite{aircraft_sos},~\cite{seiler2011susceptibility_FA18_sos},~\cite{seiler2012assessment_of_flight_sos}, among many others.

Despite the enormous impact of sum of squares optimization on polynomial optimization and related areas, the applicability of this methodology has always been limited by a single fundamental challenge, which is \emph{scalability}. Indeed, when $p(x)$ has $n$ variables and is of degree $2d$, the size of $Q$ in (\ref{eq:SOS.SDP}) is ${n+d\choose d} \times {n+d\choose d}$. Poor scaling with problem dimension is not the only difficulty here: even when their size is not too large, SDPs are arguably the most expensive class of convex optimization problems to solve. This has led the optimization community to sometimes perceive semidefinite programming as a powerful theoretical tool, but not a practical one.

\textbf{Outline.} In this paper, we review two new techniques that aim to make SOS programs more scalable. In Section \ref{sec:chordal.sparsity}, we present two efficient first-order methods based on the alternating direction method of multipliers (ADMM) to solve SDPs arising from sum of squares programming. Both methods exploit sparsity to increase the computational efficiency of SOS programs: the first method exploits the inherent sparsity of the SDPs obtained from SOS programs and can be used for any SOS program; the second one requires additional problem structure, namely that the polynomials at hand be \emph{chordally sparse}. In Section \ref{sec:dsos.sdsos}, we present techniques that do away with semidefinite programming altogether. Instead, the semidefinite program that we wish to solve is replaced by a sequence of linear or second order cone programs, which are much more tractable than SDPs. Generating these sequences amounts to constructing a series of linear and second-order cone programming-representable cones which inner approximate the set of SOS polynomials. This leads to the following conceptual question: how well can second order cone programming based techniques perform for inner approximating the set of SOS polynomials? Could we maybe even exactly represent the SOS cone using second order-representable cones? In Section \ref{sec:lower}, we review recent negative results by Fawzi which show that an exact representation is not possible, in a case as basic as the case of univariate quartics.

\textbf{Notation.} Unless otherwise specified, we will be considering throughout polynomials $p$ of degree $2d$ and in $n$ variables. We write:
\begin{equation}\label{E:Polynomial}
    p(x) = \sum_{\alpha\in \mathbb{N}^n_{2d}} p_{\alpha} x^{\alpha}, \quad p_{\alpha} \in \mathbb{R},
\end{equation}
where $x^{\alpha}=x_1^{\alpha_1}\ldots x_n^{\alpha_n}$ is a monomial of degree $|\alpha|=\sum_i \alpha_i$ and $\mathbb{N}^n_{2d}=\{\alpha \in \mathbb{N}^n~|~|\alpha|\leq 2d\}.$ We denote by $PSD_{n,2d}$ (resp. $SOS_{n,2d}$) the set of nonnegative (resp. sum of squares polynomials) in $n$ variables and of degree $2d$.
{\color{black}We further denote by $\mathbb{S}^k$ the cone of $k \times k$ symmetric matrices, and by $\mathbb{S}^{k}_+$ the cone of $k \times k$ positive semidefinite matrices.}

\section{Exploiting Sparsity in SOS Programs} \label{sec:chordal.sparsity}

In this section, we introduce two strategies that exploit sparsity to increase the computational efficiency of SOS programs. The first strategy exploits sparsity in the coefficient matching conditions arising from SOS programs for general polynomials, and the second one takes advantage of chordal sparsity for sparse polynomials. Both of them use a first-order operator splitting algorithm, known as the alternating direction method of multipliers (ADMM)~\cite{boyd2011distributed}, to efficiently compute a solution of the SDP from SOS programs at the cost of reduced accuracy.

\subsection{Sparsity in the coefficient matching conditions}

Consider a real polynomial $p(x)$ of degree $2d$ in~\eqref{E:Polynomial}. %
As mentioned in (\ref{eq:charac.sos.sdp}), if $p(x)$ is SOS, then we have
\begin{equation} \label{E:SOSdecomposition}
p(x) = \sum_{i=1}^r f_i^2(x) = \sum_{i=1}^r \left(q_i^Tz(x)\right)^2 = z^T(x)Qz(x),
\end{equation}
where $Q \succeq 0$, and $z(x)$ is a monomial basis. Generally, $z(x)$ is the vector of all monomials of degree no greater than $d$:
\begin{equation} \label{E:MonomialBasis}
z(x) = [1, x_1, x_2, \ldots, x_n, x_1^2, x_1x_2, \ldots,x_n^d]^T.
\end{equation}

Let $A_\alpha $ be the indicator matrix for the monomials $x^{\alpha}$ in the rank-one matrix $z(x)z(x)^T$. The SOS constraint \eqref{E:SOSdecomposition} can then be reformulated as
$$p(x) = \langle z(x)z(x)^T, Q \rangle
= \sum_{\alpha\in\mathbb{N}_{2d}^n} \langle A_\alpha, Q\rangle x^\alpha.$$
Matching the coefficients of the left- and right-hand sides gives the equality constraints
\begin{equation} \label{E:CoeffCondition1}
\langle A_\alpha, Q\rangle = p_{\alpha} \quad \forall \: \alpha \in \mathbb{N}^n_{2d}.
\end{equation}
These equalities~\eqref{E:CoeffCondition1} are referred to as \textit{coefficient matching conditions}~\cite{zheng2017exploiting}. Then, the existence of an SOS decomposition for $p(x)$ can be checked by solving the feasibility SDP~\cite{parrilo2003semidefinite}
\begin{equation}
\label{E:SDPSOS}
\begin{aligned}
\text{find}\quad &Q \\
\text{subject to} \quad  & \langle A_{\alpha}, Q \rangle = p_{\alpha},
\quad \alpha\in\mathbb{N}_{2d}^n,\\
& Q \succeq 0.
\end{aligned}
\end{equation}
As mentioned in Section \ref{sec:intro}, the dimension of the positive semidefinite variable $Q$ in~\eqref{E:SDPSOS} is $ \begin{pmatrix}\!\begin{smallmatrix}  n+d \\ d \\ \end{smallmatrix} \end{pmatrix} \times  \begin{pmatrix}\!\begin{smallmatrix}  n+d \\ d \\ \end{smallmatrix} \end{pmatrix} $, which grows quickly as $n$ or $d$ increases. Note that this number may be reduced by taking advantage of the structural properties of $p(x)$ to eliminate redundant monomials in $z(x)$; well-known techniques include Newton polytope~\cite{reznick1978extremal}, diagonal inconsistency~\cite{lofberg2009pre}, and symmetry property~\cite{gatermann2004symmetry}. Also, the size of the underlying SDP was investigated for some classes of matrix polynomials with sparsity in~\cite{chesi2017complexity}.

\begin{table*}
	\centering
	\setlength{\abovecaptionskip}{0pt}
	\setlength{\belowcaptionskip}{1em}
	\renewcommand\arraystretch{0.9}
	\caption{Density of nonzero elements in the equality constraints of SDP~\eqref{E:SDPSOS}}
	\label{T:density}
	\begin{tabular}{c| c c c c c c c}
		\hline \toprule[1pt] 
		$n$   & 4 & 6  & 8 & 10 & 12 & 14 & 16\\
		\hline\\[-0.75em]
		$2d = 4$ & $1.42\times 10^{-2}$ & $4.76\times 10^{-3}$ & $2.02\times 10^{-3}$ & $9.99\times 10^{-4}$  & $5.49\times 10^{-4}$& $3.27\times 10^{-4}$ & $2.06\times 10^{-4}$\\
		$2d = 6$ & $4.76\times 10^{-3}$ & $1.08\times 10^{-3}$ & $3.33\times 10^{-4}$ & $1.25\times 10^{-4}$  & $5.39\times 10^{-5}$& $2.58\times 10^{-5}$ & $1.34\times 10^{-5}$\\
		$2d = 8$ & $2.02\times 10^{-3}$ & $3.33\times 10^{-4}$ & $7.77\times 10^{-5}$ & $2.29\times 10^{-5}$  & $7.94\times 10^{-6}$ & $3.13\times 10^{-6}$ & $1.36\times 10^{-6}$\\
		\bottomrule[1pt]
	\end{tabular}
\end{table*}

One important feature of~\eqref{E:SDPSOS} is that the coefficient matching conditions are sparse, in the sense that each equality constraint in~\eqref{E:SDPSOS} only involves a small subset of entries of $Q$ \cite{zheng2017exploiting}, since only a small subset of entries of the product $z(x)z(x)^T$ are equal to a given monomial $x^\alpha$. In particular, we re-index the constraint matching conditions~\eqref{E:CoeffCondition1} using integer indices $i=1,\ldots,m$, where $m = \begin{pmatrix}\!\begin{smallmatrix}   n+2d \\ 2d \\ \end{smallmatrix} \end{pmatrix}$. Let $\mathrm{vec}: \mathbb{S}^N \to \mathbb{R}^{N^2}$ be the operator mapping a matrix to the stack of its columns, and define
$$
    A = \begin{bmatrix} \text{vec}(A_1), \cdots, \text{vec}(A_m) \end{bmatrix}^T.
$$
Then, the equality constraints in~\eqref{E:SDPSOS} can be rewritten as
\begin{equation} \label{E:SDPSOSequaltiy}
    A\cdot \text{vec}(Q) = b,
\end{equation}
where $b\in \mathbb{R}^m$ is a vector collecting the coefficients $p_{i}$ of $p$. We have the following result.
\begin{theorem}[Sparsity of constraints~\cite{zheng2017exploiting}] \label{T:sparsity}
	Let $A$ be the coefficient matrix of the equality constraints for~\eqref{E:SDPSOS}. The density of nonzero elements in $A$ is $\mathcal{O}({\frac{1}{n^{2d}}})$.
\end{theorem}

Note that this result holds for dense polynomials. 
Also, the density decreases quickly as  $n$ or $d$ increases, which means the SDP~\eqref{E:SDPSOS} becomes very sparse for large-scale (dense or not) polynomials. In Table~\ref{T:density}, we list the density of nonzero elements in typical cases. Therefore, it is desirable to exploit this sparsity to improve the computational efficiency of~\eqref{E:SDPSOS}. Let us represent $A = [a_1, a_2, \ldots, a_m]^T$, so that each vector $a_i$ is a row of $A$, and let $H_i, i=1,\ldots,m$ be ``entry-selector'' matrices of 1's and 0's that select the nonzero elements of $a_i$. We have the following equivalence.
\begin{equation} \label{E:LocalVariable}
A\cdot \text{vec}(Q) = b \Leftrightarrow
\begin{cases}
(H_ia_i)^Tz_i = b_i  \\
z_i = H_i\cdot \text{vec}(Q)
\end{cases},i = 1, \ldots, m,
\end{equation}
where $z_i$ is a copy of the non-zero elements of $\text{vec}(Q)$ in the $i$-th equality constraint. Note that only the non-zero elements are involved in~\eqref{E:LocalVariable}. Also, the equalities in~\eqref{E:LocalVariable} are enforced \emph{individually} by $z_i$, not \emph{simultaneously} as in~\eqref{E:SDPSOSequaltiy}.

Together with equivalence~\eqref{E:LocalVariable} that exploits the sparsity in $A$, we can apply ADMM to~\eqref{E:SDPSOS}, resulting in an efficient algorithm that is free of any matrix inversion. Each iteration of the resulting ADMM algorithm consists of one conic projection and multiple quadratic programs with closed-form solutions; we refer the interested reader to~\cite{zheng2017exploiting} for details.  

\subsection{Chordal sparsity in sparse polynomials}

The strategy that exploits the sparsity in the coefficient matching conditions works for general (dense or not) polynomials. However, the dimension of $Q$ in~\eqref{E:SDPSOS} is unchanged, which may still require extensive computation for large-scale instances. If the polynomial $p(x)$ has structured sparsity, the large cone constraint $Q\succeq0$ can be replaced by a set of smaller cone constraints~\cite{fukuda2001exploiting,waki2006sums,zheng2016fast}. Specifically, we assume the polynomial $p(x)$ has \emph{correlative sparsity}, first introduced by Waki \emph{et al.}~\cite{waki2006sums}, which is a higher level of view of the sparsity of a polynomial in terms of the interaction of variables $x_i$. Given a polynomial $p(x)$ in~\eqref{E:Polynomial}, the correlative sparsity pattern is represented by a matrix $R \in \mathbb{S}^n$, 
\begin{equation*}
R_{ij} = \begin{cases}
1, \quad \text{if } i = j \text{ or } \alpha_i, \alpha_j \geq 1, p_{\alpha} \neq 0 \\
0, \quad \text{otherwise} \\\end{cases}.
\end{equation*}
Further, we can associate an undirected graph $\mathcal{G}(\mathcal{V},\mathcal{E})$ with $\mathcal{V} = \{1, \ldots,n\}$ and $\mathcal{E} =\{(i,j) \mid R_{ij} = 1, i \leq j\}$. Then, Waki \emph{et al.} proposed that multiple sets of monomial basis could be used to construct an SOS polynomial~\cite{waki2006sums}, \emph{i.e.},
\begin{equation} \label{E:SOSdecompositionChordal}
p(x) = \sum_{i=1}^q z_i^T(x)Q_iz_i(x),
\end{equation}
where $z_i(x)$ is a monomial basis and $Q_i \succeq 0$. The choice of $z_i(x)$ depends on the maximal clique (see the precise definition below) of a \emph{chordal extension} of the graph $\mathcal{G}(\mathcal{V},\mathcal{E})$, and the dimension of $Q_i$ becomes small if the size of the largest maximal clique is small (please refer to~\cite{waki2006sums} for further information).

As we have seen that an SOS program can always be transformed into a special SDP (see Section \ref{sec:intro}), we focus in the following on a general sparse SDP and introduce an ADMM algorithm that exploits the inherent chordal sparsity. We consider the following primal standard SDP.
\begin{equation}
\label{E:PrimalSDP}
\begin{aligned}
\min_{X} \quad & \langle C, X \rangle \\
\text{subject to} \quad & \langle A_i,X \rangle = b_i, i = 1, \ldots, m,\\
& X\succeq 0.
\end{aligned}
\end{equation}
For the sake of completeness, we first introduce several graph-theoretic concepts. Given an undirected graph $\mathcal{G}(\mathcal{V},\mathcal{E})$, a subset of vertices $\mathcal{C} \subseteq \mathcal{V}$ is called a clique if $(i,j) \in \mathcal{E}, \forall \: i,j \in \mathcal{C}, i \neq j$. The clique is called maximal if it is not a subset of any other clique. An undirected graph $\mathcal{G}$ is called \emph{chordal} if every cycle of length greater than three has at least one chord.  Note that if $\mathcal{G}(\mathcal{V}, \mathcal{E})$ is not chordal, it can be \emph{chordal extended}, \emph{i.e.}, we can construct a chordal graph $\mathcal{G}'(\mathcal{V}, \mathcal{E}') $ by adding additional edges to $ \mathcal{E} $. More details can be found in~\cite{blair1993introduction}.

Let $\mathcal{G}(\mathcal{V},\mathcal{E})$ be an undirected graph with self-loops. We say that $X$ is a partial symmetric matrix defined by $\mathcal{G}$ if $X_{ij} = X_{ji}$ are given when $(i,j) \in \mathcal{E}$, and arbitrary otherwise. We define the following sparse cones.
\begin{equation*}
\begin{aligned}
\mathbb{S}^n(\mathcal{E},?) = &\{ X \in \mathbb{S}^n \mid X_{ij} =X_{ji}\text{ given if } (i,j) \in \mathcal{E}  \}, \\
\mathbb{S}_{+}^n(\mathcal{E},?) = &\{ X \in \mathbb{S}^n(\mathcal{E},?) \mid
\notag \\ &\qquad\,\,\quad\exists M \succeq 0, \, M_{ij} = X_{ij} , \forall (i,j) \in \mathcal{E}  \}.
\end{aligned}
\end{equation*}
Given a clique $\mathcal{C}_k$ of $\mathcal{G}$, we let $E_{\mathcal{C}_k} \in \mathbb{R}^{\mid \mathcal{C}_k\mid \times n}$ be the matrix with $(E_{\mathcal{C}_k})_{ij} = 1$ if $\mathcal{C}_k(i) = j$ and zero otherwise, where $\mathcal{C}_k(i)$  is the $i$-th vertex in $\mathcal{C}_k$, sorted in the natural ordering. Then, we have the following result.
\begin{theorem} [Grone's theorem~\cite{grone1984positive}]\label{T:ChordalCompletionTheorem}
	Let $\mathcal{G}(\mathcal{V},\mathcal{E})$ be a chordal graph with a set of maximal cliques $\{\mathcal{C}_1,\mathcal{C}_2, \ldots, \mathcal{C}_p\}$. Then, $X\in\mathbb{S}^n_+(\mathcal{E},?)$ if and only if
	$X_k = E_{\mathcal{C}_k} X E_{\mathcal{C}_k}^T \in \mathbb{S}^{\vert \mathcal{C}_k \vert}_+$
	for all $k=1,\,\ldots,\,p$.
\end{theorem}
\begin{remark}
	This theorem allows us to equivalently replace $\mathbb{S}^n_{+}(\mathcal{E},?)$ with a set of coupled but smaller convex cones. A dual result can be found in~\cite{agler1988positive}. These results have been exploited in interior-point methods for SDPs~\cite{fukuda2001exploiting}; also, see recent applications in stability analysis and controller synthesis of large-scale linear systems~\cite{mason2014chordal, zheng2016scalable, zheng2016chordal}.
\end{remark}

We assume that~\eqref{E:PrimalSDP} is sparse with an \emph{aggregate sparsity pattern} described by $\mathcal{G}(\mathcal{V},\mathcal{E})$, meaning that $(i,j)\in\mathcal{E}$ if and only if the entry $ij$ of at least one of the matrices $C,\,A_1,\,\ldots,\,A_m$ is nonzero. Also, it is assumed that $\mathcal{G}$ is chordal with a set of maximal cliques $\mathcal{C}_1,\ldots,\mathcal{C}_p$. In~\eqref{E:PrimalSDP}, only the entries of $X$ corresponding to the edges $\mathcal{E}$ appear in the cost and constraint functions. Therefore, the constraint $X\in\mathbb{S}^n_+$ can be replaced by $X\in\mathbb{S}^n_+(\mathcal{E},?)$. Using Theorem~\ref{T:ChordalCompletionTheorem}, we can reformulate~\eqref{E:PrimalSDP} as
\begin{equation}
\label{E:DecomposedPrimalSDP}
\begin{aligned}
\min_{X,X_1,\ldots,X_p} \quad & \langle C, X \rangle \\
\text{subject to} \quad  &\langle A_i,X \rangle = b_i, && i = 1, \ldots,m\\
& X_k - E_k X E_k^T = 0, &&k=1,\,\ldots,\,p,\\
& X_k \in \mathbb{S}^{\vert \mathcal{C}_k\vert}_{+}, &&k=1,\,\ldots,\,p.
\end{aligned}
\end{equation}
In other words, the original large semidefinite cone is decomposed into multiple smaller cones at the cost of introducing a set of consensus constraints between the variables. Together with this reformulation that exploits the aggregate sparsity pattern to reduce the cone dimension, we can apply ADMM to~\eqref{E:DecomposedPrimalSDP}, which results in an efficient algorithm that works with smaller positive semidefinite cones; see~\cite{zheng2016fast} for details. Similar decompositions are available for the dual standard SDP and the homogeneous self-dual embedding of sparse SDPs~\cite{ZFPGW2017chordal,zheng2016hsde}.

\subsection{Numerical results}

The two strategies have been implemented in the MATLAB packages SOSADMM and CDCS~\cite{CDCS}, respectively. These two packages are available from
\urlstyle{tt} \small \url{https://github.com/oxfordcontrol/SOSADMM} \normalsize
and
\urlstyle{tt} \small \url{https://github.com/oxfordcontrol/CDCS}\normalsize.
This section presents numerical tests of SOSADMM and CDCS on the random unconstrained polynomial optimization problems (more numerical results can refer to~\cite{ZFPGW2017chordal,zheng2017exploiting}). In the experiments, we set the termination tolerance to $10^{-4}$, and the maximum number of iterations to $2\times 10^3$ for SOSADMM and CDCS. The primal method in CDCS was used, and SeDuMi~\cite{sturm1999using} was used as a benchmark solver. Consider the polynomial minimization problem
$ \min_{x\in\mathbb{R}^n}  p(x)$, where $p(x)$ is a given polynomial. As described in Section \ref{sec:intro}, we can obtain an SDP relaxation as
\begin{equation}
\label{E:POPSDP}
\begin{aligned}
\max \quad & \gamma \\
\text{subject to} \quad  &p(x) - \gamma \quad \text{is SOS}.
\end{aligned}
\end{equation}
We generated $p(x)$ according to $ p(x) = p_0(x) + \sum_{i=1}^n x_i^{2d}$, where $p_0(x)$ is a random polynomial with normally distributed coefficients of degree strictly less than $2d$. We used GloptiPoly~\cite{henrion2003gloptipoly} to generate the examples. Table~\ref{T:TimePOP} compares the CPU time (in seconds) required to solve the SOS relaxation as the number of variables was increased $n$ with $d=2$. Both SOSADMM and CDCS-primal were faster than SeDuMi on these examples. Also, the optimal value returned by SOSADMM was within 0.05\% of the high-accuracy value returned by SeDuMi. Note that SOSADMM was faster than CDCS-primal in the experiments and this is expected since the random polynomials were dense; major computational improvements have been achieved by SOSADMM. For brevity, the interested reader is referred to~\cite{ZFPGW2017chordal} for more numerical results of CDCS on sparse SDPs.%

\begin{table}[t]
	\centering
	\setlength{\abovecaptionskip}{0pt}
	\setlength{\belowcaptionskip}{0em}
	\renewcommand\arraystretch{1.1}
	\caption{CPU time (s) to solve the SDP relaxations \eqref{E:POPSDP}. $N$ is the size of the PSD cone, $m$ is the number of constraints.}
	\label{T:TimePOP}
	\begin{tabular}{c c c|  c c c }
		\hline \toprule[1pt] 
		\multicolumn{3}{c|}{Dimensions} &  \multicolumn{3}{c}{CPU time (s)} \\
		\hline
		$n$   & $N$ & $m$   & SeDuMi &
		\begin{tabular}[x]{@{}c@{}}CDCS\\[-0.25em](primal)\end{tabular} &
		\begin{tabular}[x]{@{}c@{}}SOS-\\[-0.25em]ADMM\end{tabular} \\
		$2$ & 6 & 14  & 0.108  & 0.163 & 0.041  \\
		$6$ & 28 & 209  & 0.295  & 0.212 & 0.093  \\
		$10$ & 66 & 1000  & 4.197  & 0.340 & 0.294  \\
		$14$ & 120 & 3059  & 53.68   & 0.575 & 0.490 \\
		$18$ & 190 & 7314  & 621.2 & 1.696  & 1.339  \\
		$20$ & 231 & 10625  & 1806.6  & 4.694 & 2.362  \\
		\bottomrule[1pt]
	\end{tabular}
\end{table}

\section{Linear and second-order programming-based alternatives to SOS programs} \label{sec:dsos.sdsos}

In this section, we focus on another class of methods that enables us to increase the computational efficiency of SOS programs. These consist in replacing the semidefinite program that underlies any SOS program (see Section \ref{sec:intro}) by more tractable convex programs such as linear or second-order cone programs. 

\subsection{DSOS and SDSOS programs}\label{sec:dsos.sdsos.intro}

Recall from Section \ref{sec:intro} that a polynomial $p(x)$ of degree $2d$ and in $n$ variables is SOS if and only if there exists a positive semidefinite matrix $Q$ such that $p(x)=z(x)^TQz(x)$, where $z(x)=[1,x_1,x_2,\ldots,x_n^d]$ is the vector of monomials of degree $d$. As a consequence, solving any SOS program amounts to solving a semidefinite program where the variable {\color{black}$Q$ is of size $\binom{n+d}{d}$}. As the number of variables and degree of the polynomial $p(x)$ at hand increase, the cost of solving such a semidefinite program can quickly become prohibitive. The idea proposed in~\cite{iSOS_journal}  is to replace the condition that the Gram matrix $Q$ be positive semidefinite with stronger but cheaper conditions in the hope of obtaining more efficient inner approximations to the cone $SOS_{n,2d}$. Two such conditions come from the concepts of \emph{diagonally dominant} and \emph{scaled diagonally dominant} matrices in linear algebra. We recall these definitions below.

\begin{definition}\label{def:dd.sdd}
	A symmetric matrix $A$ is \emph{diagonally dominant} (dd) if $a_{ii} \geq \sum_{j \neq i} |a_{ij}|$ for all $i$. We say that $A$ is \emph{scaled diagonally dominant} (sdd) if there exists a diagonal matrix $D$, with positive diagonal entries, which makes $DAD$ diagonally dominant.
\end{definition}

We refer to the set of $n \times n$ dd (resp. sdd) matrices as $DD_n$ (resp. $SDD_n$). The following inclusions are a consequence of Gershgorin's circle theorem \cite{gersh}:
$$DD_n\subseteq SDD_n\subseteq \mathbb{S}_+^n.$$

\begin{definition}[\cite{iSOS_journal}] \label{def:dsos.sdsos.rdsos.rsdsos}
	A polynomial $p(x)$ of degree $2d$ is said to be \emph{diagonally-dominant-sum-of-squares} (DSOS) if it admits a representation as {\color{black} $p(x)=z^T(x)Qz(x)$,} where $Q$ is a dd matrix and $z(x)$ is the standard vector of monomials of degree $d$. A polynomial $p(x)$ of degree $2d$ is said to be \emph{scaled-diagonally-dominant-sum-of-squares} (SDSOS) if it admits a representation as {\color{black}$p(x)=z^T(x)Qz(x)$}, where $Q$ is an sdd matrix and $z(x)$ is the standard vector of monomials of degree $d$.
	
\end{definition}

\begin{figure}[h]
	\centering
	\includegraphics[width=0.4\textwidth]{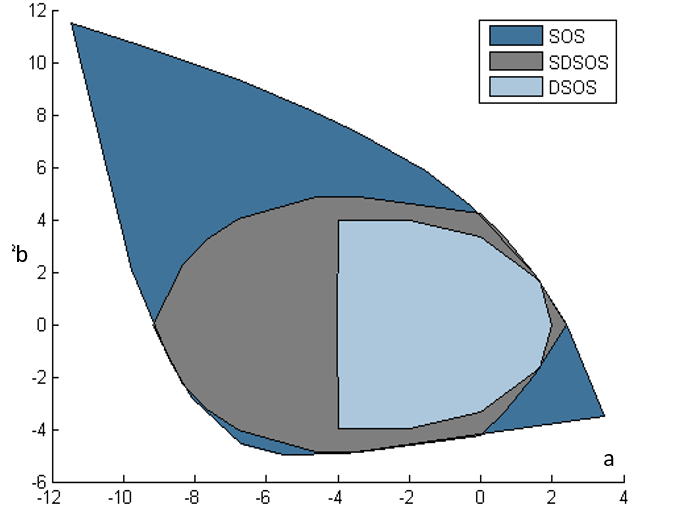}
	\caption{A comparison of the DSOS/SDSOS/SOS cones on an example}
	\label{fig:s.d.sos}
\end{figure}

Let us denote the cone of polynomials in $n$ variables and degree $2d$ that are DSOS and SDSOS by $DSOS_{n,2d}$, $SDSOS_{n,2d}$. The following inclusion relations are straightforward: $$DSOS_{n,2d}\subseteq SDSOS_{n,d}\subseteq SOS_{n,2d}\subseteq PSD_{n,2d}.$$ 
An illustration of how sections of these different cones compare on an example is given in Figure~\ref{fig:s.d.sos}, taken from \cite{iSOS_journal}. 
We consider a parametric family of polynomials parameterized by $a$ and $b$, $$p_{a,b}(x_1,x_2)=2x_1^4+2x_2^4+ax_1^3x_2+(1-a)x_1^2 x_2^2+bx_1x_2^3,$$ and plot in the figure the values of $(a,b)$ for which the polynomial is DSOS (innermost set), SDSOS (set containing the DSOS set), and SOS (or equivalentally nonnegative in this case) (outermost set).

These definitions give rise to the notion of DSOS and SDOS programs. With similar notation to (\ref{eq:SOS.programs}), we say that a DSOS (resp. SDSOS) program is an optimization problem of the following form:
\begin{equation}\label{eq:DSOS.SDSOS.programs}
\begin{aligned}
&\min_{p \in \mathbb{R}_{n,2d}[x]} &&C(p)\\
&\text{s.t. } &&A(p)=b\\
& &&p \text{ is DSOS (resp. SDSOS)}
\end{aligned}
\end{equation}
Note that (\ref{eq:DSOS.SDSOS.programs}) provides upperbounds on the optimization problem given in (\ref{eq:SOS.programs}) as the set of DSOS and SDSOS polynomials is a subset of the set of SOS polynomials. This loss in solution accuracy is compensated by gains in terms of scalability and solving-time. This is a consequence of the following theorem.

\begin{theorem}[\cite{iSOS_journal}] For any fixed $d$, solving a DSOS (resp. SDSOS) program can be done with linear programming (resp. second order cone programming) of size polynomial in $n$.
\end{theorem}

The ``LP part'' of this theorem is not hard to see. The equality $p(x)=z(x)^TQz(x)$ gives rise to linear equality constraints between the coefficients of $p$ and the entries of the matrix $Q$. The requirement of diagonal dominance on the matrix $Q$ can also be described by linear inequality constraints on $Q$. The ``SOCP part'' of the statement is not as straightforward and its proof can be found in ~\cite{iSOS_journal}.

We illustrate the gains that one can make using these methods in Table \ref{tab:dsos.sdsos} taken from \cite{iSOS_journal}. We have reported the time and bounds obtained when mininimizing a random polynomial of degree $d=4$ and with a varying number of variables $n$ over the unit sphere, using a DSOS, SDSOS and SOS program (see~\cite{iSOS_journal} for the precise formulations). Note that when $n$ is small, the SOS program returns a better bound in slightly longer times than the DSOS/SDSOS programs. However, when $n$ gets large, the SOS program cannot be solved due to memory issues whereas both the DSOS and SDSOS programs run in the order of seconds. These results were obtained on a 3.4 GHz Windows computer with 16 GB of memory.

\begin{table*}[t]
	5\small
    \setlength{\abovecaptionskip}{0pt}
	\setlength{\belowcaptionskip}{0em}
	\renewcommand\arraystretch{1.1}
    \caption{ Lower bounds obtained using S/D/SOS programs to compute the minimum of a quartic form on the sphere for varying $n$, along with run times (in secs).}
		\label{tab:dsos.sdsos}
	\begin{center}
		\begin{tabular}{ c  c c  c  c c c  c c c  c c c c c}
			\hline \toprule[1pt]
			& \multicolumn{2}{c}{$n=15$}  & & \multicolumn{2}{c}{$n=20$} & &\multicolumn{2}{c}{$n=25$} & & \multicolumn{2}{c}{$n=30$} & & \multicolumn{2}{c}{$n=40$}\\
			\cline{2-3} \cline{5-6} \cline{8-9} \cline{11-12} \cline{14-15}
			& bd & t(s) && bd & t(s) && bd & t(s) && bd & t(s) && bd & t(s) \\
			\cline{2-3} \cline{5-6} \cline{8-9} \cline{11-12} \cline{14-15}
			DSOS & -10.96 & 0.38 & & -18.012 & 0.74 && -26.45 & 15.51 &&  -36.85 & 7.88 &&  -62.30 & 10.68  \\ 
			SDSOS  &  -10.43 & 0.53 &&  -17.33 & 1.06 &&  -25.79 & 8.72& & -36.04 & 5.65 &&  -61.25 & 18.66   \\ 
			SOS & -3.26 & 5.60 && -3.58 & 82.22 & & -3.71 & 1068.66 && NA & NA && NA & NA   \\
            \bottomrule[1pt]
		\end{tabular}
	\end{center}
\end{table*}

\subsection{Improving on DSOS and SDSOS programming}

As mentioned previously, the advantages of substituting an SOS program with a DSOS or SDSOS program are scalability and computational efficiency of the program obtained. This comes at the cost of solution accuracy. In this section, we present two methods for mitigating the loss in accuracy that we observe. These methods involve constructing a sequence of iterative linear or second-order cone programs where the first iteration consists in solving the DSOS/SDSOS program given in (\ref{eq:DSOS.SDSOS.programs}). 
Our goal throughout will be to solve the optimization problem given in (\ref{eq:SOS.programs}).

\subsubsection{Column generation method \cite{AAASDGH}} \label{sec:col.gen} For simplicity, we present here the linear programming-based version of the algorithm. An analogous method based on second-order cone programming can be found in \cite{AAASDGH}. To understand this method, the following characterization of diagonally dominant matrices is needed \cite{dd_extreme_rays}: A symmetric matrix $M$ is diagonally dominant if and only if it can be written as
{\color{black} 
$$
    M=\sum_{i=1}^{n^2} \alpha_i v_iv_i^T, \alpha_i\geq 0,
$$}where $\{v_i\}$ is the set of all nonzero vectors in $\mathbb{R}^n$ with at most $2$ nonzero components, each equal to $\pm 1$. The first iteration of our algorithm (i.e., solving (\ref{eq:DSOS.SDSOS.programs})) then amounts to solving:
\begin{equation} \label{eq:col.gen}
\begin{aligned}
&\min_{p\in \mathbb{R}_{n,2d}[x],\alpha} &&C(p)\\
&\text{s.t. } &&A(p)=b\\
& &&p(x)=z(x)^TQz(x), \forall x,\\
& &&Q=\sum_{i} \alpha_iv_iv_i^T, \alpha_i \geq 0,
\end{aligned}
\end{equation}
where $\{v_i\}$ are fixed. At each iteration, one adds a new ``column" $v$ to the set $\{v_i\}$ and the problem is then solved again. This leads to Algorithm \ref{alg:alg.col.gen}.

\begin{algorithm}
	\caption{Column generation algorithm}\label{alg:alg.col.gen}
	\begin{algorithmic}[1]
		\State  \textbf{initialize:} Solve (\ref{eq:col.gen}). Obtain $\alpha$ and $p$.
		\Repeat\\
		Find a vector $v$ using the dual as described below \\
		Solve \newline $$\begin{matrix}\min_{\alpha,p} C(p)\\ \text{s.t. } A(p)=b\\ p(x)=z(x)^TQz(x), \forall x,\\ Q=\sum_{i}\alpha_i v_iv_i^T+\alpha vv^T \\ \alpha_i \geq 0~\forall i, \alpha \geq 0 \end{matrix} $$
		\Until Termination Condition is met\\
		\Return $C(p)$ and $p$
	\end{algorithmic}
\end{algorithm}

Note that at each iteration, this algorithm can only improve: indeed, by taking $\alpha=0$, one recovers the solution to the previous iteration. To obtain strict improvement, one needs to carefully pick the vector $v$ that we add to our set of columns. One way of doing this is via the dual of the problem given in (\ref{eq:col.gen}). The constraint $Q=\sum_{i} \alpha_i v_iv_i^T$ in the primal gives rise to a constraint of the type $v_i^TXv_i \geq 0, ~\forall i$ in the dual, where $X$ is a dual variable. A good choice of a vector $v$ is then given by any vector $v$ such that $v^TXv<0$, see \cite{AAASDGH} for more information regarding the choice of $v$.

One can choose to terminate the algorithm under different conditions, e.g., lack of improvement of the optimal value/solution, or expiration of the time/computational budget associated to the task of solving the SDP.

An illustration of the performance of this technique is given in Table \ref{tab:min_triples} taken from \cite{AAASDGH}. The setting is analogous to the one used to obtain Table \ref{tab:dsos.sdsos}: we report the time and bounds obtained when minimizing a random degree $d=4$ homogeneous polynomial over the unit sphere. The experiments were run on a 2.33 GHz Linux machine with 32 GB of memory. In each iteration, we add an appropriate new vector $v$ to the sequence $\{v_i\}$ that has at most three nonzero elements, each equaling $1$ or $-1$. We stop the algorithm when either of the two following conditions are met: lack of improvement in the optimal value or time budget of 600s exceeded. Note that one can significantly improve on the initial approximations provided in Table \ref{tab:dsos.sdsos} within a reasonable 10~min time lapse.

%
%

\begin{table*}
    \setlength{\abovecaptionskip}{0pt}
	\setlength{\belowcaptionskip}{0em}
	\renewcommand\arraystretch{1.2}
    \caption{ Lower bounds obtained using column generation to compute the minimum of a quartic form on the sphere for varying $n$, along with run times (in secs).}
	\label{tab:min_triples}
	\begin{center}
		\begin{tabular}{  c  c c  cccc c c  c  c c  c c c c }
			\hline \toprule[1pt]
			 &  \multicolumn{2}{c}{$n=15$}  &&  \multicolumn{2}{c}{$n=20$} && \multicolumn{2}{c}{$n=25$} &&  \multicolumn{2}{c}{$n=30$} && \multicolumn{2}{c}{$n=40$} \\
			\cline{2-3} \cline{5-6} \cline{8-9} \cline{11-12} \cline{14-15}
			\ & bd & t(s) && bd & t(s) && bd & t(s) && bd & t(s) && bd & t(s) \\
			\cline{2-3} \cline{5-6} \cline{8-9} \cline{11-12} \cline{14-15}			
			Col Gen & $-5.57$ & $31.19$ &&$-9.02$ & $471.39$ && $-20.08$ & $600$ &&$-32.28$ & $600$ && $-35.144$ & $600$  \\
    \bottomrule[1pt]
		\end{tabular}
	\end{center}
\end{table*}

\subsubsection{Sum of squares basis pursuit \cite{AAA_GH_Basis_Pursuit}}

The idea behind this algorithm is the following. Let us assume that one could solve the problem given in (\ref{eq:SOS.programs}) and obtain the Gram matrix $Q^*$ associated to the optimal polynomial $p^*$. If we changed the monomial basis $z(x)$ in (\ref{eq:SOS.SDP}) to the monomial basis $L^*z(x)$, where $L^*$ is the Cholesky decomposition (or square root) of $Q^*$, then we could use linear programming to recover the optimal solution of (\ref{eq:SOS.programs}). Indeed, in this case, the optimal Gram matrix would be given by the identity matrix, which is diagonal, and searching for a diagonal Gram matrix can be done via linear programming. In our case, we do not have access to $L^*$, but the idea is to work towards such a basis. This procedure is detailed in Algorithm \ref{alg:bas.pursuit}.

\begin{algorithm}
	\caption{Sum of squares basis pursuit algorithm}\label{alg:bas.pursuit}
	\begin{algorithmic}[1]
		\State  \textbf{initialize:} Solve (\ref{eq:DSOS.SDSOS.programs}). Obtain the Gram matrix $Q$ associated to the optimal $p$. Compute $L=chol(Q)$.
		\Repeat\\
		Solve \newline $$\begin{matrix}\min_{\alpha,p} C(p)\\ \text{s.t. } A(p)=b\\ p(x)=z(x)^TL^TQLz(x), \forall x,\\ Q \text{ is d.d./s.d.d.}  \end{matrix} $$\\
		Obtain the Gram matrix $Q$ associated to the optimal $p$. Set $L \leftarrow chol(Q)L.$
		\Until Termination Condition is met\\
		\Return $C(p)$ and $p$
	\end{algorithmic}
\end{algorithm}

The algorithm will terminate under similar conditions to the ones given in Section \ref{sec:col.gen}. Note that this algorithm is guaranteed to converge: at each iteration, the optimal value of the problem decreases (indeed, setting $Q=I$ enables us to recover the solution at the previous iteration), and is lowerbounded by the optimal solution to the SDP in (\ref{eq:SOS.programs}).


We finally point the reader to very recent work~\cite{pop_hierarchy}, which goes beyond DSOS and SDSOS optimization and produces a converging hierarchy for the polynomial optimization problem that does not even require the use of linear or second order cone programming. This hierarchy only involves multiplying certain polynomials together and checking whether the coefficients of the product are nonnegative.

\section{Limitations on representing SOS cones with bounded size PSD blocks}
\label{sec:lower}
In Section~\ref{sec:dsos.sdsos.intro} we discussed methods to certify non-negativity of
polynomials by showing that they are SDSOS. Checking that a polynomial is SDSOS
involves solving a second-order cone program. It is well-known that any second-order cone program can be written as a semidefinite
program in which all blocks have size $2\times 2$ and hence testing whether a polynomial is SDSOS simply amounts to solving an SDP involving only $2\times 2$ blocks. 
Restricting to SDPs with only small blocks is attractive because such problems can be solved more
efficiently than general SDPs of the same size.

On the one hand, in Section~\ref{sec:dsos.sdsos.intro} we saw that approximating SOS cones
with SDSOS cones (and variations on this idea) are empirically very powerful.
On the other hand, in general the SOS cone strictly contains the SDSOS cone.
Beyond the SDSOS cone, there are potentially many other ways to certify
non-negativity of a family of polynomials using SDPs with only $2\times 2$
blocks. In this section we consider recent results illustrating the
\emph{limitations} of modeling with linear matrix inequalities (LMIs) having blocks of bounded size
(particularly $2\times 2$ blocks). Concretely, we consider the following
question:
\begin{description}
	\item[Q1] Is it possible to \emph{exactly represent} SOS cones using
	LMIs that only involve $2 \times 2$ blocks, (or, more generally, blocks of
	bounded size)?
\end{description}
Closely related, and more refined, is the following approximation version of
the question:
\begin{description}
	\item[Q2] Given a positive integer $p$, how well can we
	\emph{approximate} SOS cones with convex cones that can be described using LMIs
	that involve at most $p$, $2\times 2$ blocks?
\end{description}
The focus of this section is to discuss recent results of
Fawzi~\cite{fawzi2016representing} showing that the answer to Q1 is already
negative for non-negative univariate quartics. Addressing Q2 remains a research challenge.
In presenting Fawzi's result, we
briefly describe recently developed ideas and tools for reasoning about all
possible LMI descriptions of a convex set (with fixed block sizes). These are
based on the idea of \emph{cone ranks} (introduced by Gouveia, Parrilo, and
Thomas~\cite{gouveia2013lifts}) of certain non-negative matrices associated
with the convex set.

\subsection{General PSD lifts with fixed block size}
\label{sec:js-bg}

As mentioned previously, we use the notation $\mathbb{S}_+^k$ for the cone of $k\times k$ positive
semidefinite matrices, and the notation $(\mathbb{S}_+^k)^p$ for the Cartesian
product of $p$ copies of $\mathbb{S}_+^k$.

With this notation, we can now formalize the idea of a representation of a
convex set with LMIs involving only $2\times 2$ blocks.  The following
definition is a special case of a definition due to Gouveia, Parrilo, and
Thomas~\cite{gouveia2013lifts}.
\begin{definition}
	A convex cone $C\subseteq \mathbb{R}^n$ has a \emph{proper
	$(\mathbb{S}_+^2)^p$-lift}\footnote{Clearly there is an analogous definition
		for a proper $(\mathbb{S}_+^k)^p$-lift for any fixed $k$. All the definitions
		and results in Section~\ref{sec:js-bg} extend to that case.} if there is a
	subspace $L$ of $(\mathbb{S}^2)^p$ and a linear map $\pi:
	(\mathbb{S}^2)^p\rightarrow \mathbb{R}^p$ such that
	$$C = \pi\left[(\mathbb{S}_+^2)^p \cap L\right]$$
	and $L$ meets the interior of $(\mathbb{S}_+^2)^p$.
\end{definition}
More concretely, $C$ has a \emph{$(\mathbb{S}_+^2)^p$-lift} if and only if it
can be expressed in the form
\[ C = \left\{x\in \mathbb{R}^n: \exists y\in \mathbb{R}^m\;\;\textup{s.t.}\;\;
\sum_{i=1}^{n}A_ix_i + \sum_{j=1}^{m}B_jy_j \succeq 0\right\}\]
where the $A_i$ and the $B_j$ are $2p \times 2p$ symmetric matrices that are
all block diagonal (with the same block structure) consisting of $p$ blocks,
each of size $2\times 2$.

Question Q1 can now be expressed more concisely as:
\begin{description}
	\item[Q1'] For which $(n,d)$ does there exist a finite $p$ such that
	$SOS_{n,d}$ has a $(\mathbb{S}_+^2)^p$-lift?
\end{description}
Answering questions like this in the negative, i.e., showing that lifts do
\emph{not} exist, has proven very challenging.  The work of Gouveia, Parrilo,
and Thomas~\cite{gouveia2013lifts} (building on ideas of
Yannakakis~\cite{yannakakis1991expressing} in the case of linear programming),
made a connection between the existence of lifts and a certain generalization
of the non-negative rank of an entry-wise non-negative matrix.
\begin{definition}
	If $S$ is an $a\times b$ matrix with non-negative entries, the
	$\mathbb{S}_+^2$\emph{-rank} of $S$ is the smallest $p$ such that
	$$ S_{ij} = \sum_{k=1}^{p} \langle A_{ik},B_{jk}\rangle$$
	where $A_{ik},B_{jk}\in \mathbb{S}_+^2$ for all $1\leq i\leq a$,
	$1\leq j \leq b$ and $1\leq k \leq p$.
\end{definition}
There is a connection between  $(\mathbb{S}_+^2)^p$ lifts of a convex set $C$
and the $\mathbb{S}_+^2$-rank of various non-negative matrices associated with
$C$ (so-called \emph{slack matrices} of $C$, see~\cite{gouveia2013lifts}).

The following result is a special case of~\cite[Theorem 1]{gouveia2013lifts}
expressed in the conic setting.
\begin{theorem}
	\label{thm:psdrank}
	Let $v_1,\ldots,v_b\in C$ and $\ell_1,\ldots,\ell_a\in C^* = \{\ell\in
	\mathbb{R}^n: \langle \ell,v\rangle \geq 0\;\;\textup{for all $v\in C$}\}$.  If
	$C$ has a proper $(\mathbb{S}_+^2)^p$ lift, then the non-negative matrix with
	entries
	$S_{ij} = \langle \ell_i,v_j\rangle$
	has $\mathbb{S}_+^2$-rank at most $p$.
\end{theorem}
One implication of this result, is that a lower bounds on the
$\mathbb{S}_+^2$-rank of any $S$ constructed as in Theorem~\ref{thm:psdrank}
gives a lower bound on the smallest $p$ for which $C$ has an
$(\mathbb{S}_+^2)^p$-lift. To show that $C$ has no $(\mathbb{S}_+^2)^p$-lift
for any $p$, it is enough to find a sequence of points $v_1,v_2,\ldots\in C$
and $\ell_1,\ell_2,\ldots\in C^*$, such that the corresponding sequence of
non-negative matrices $S$ (of growing dimensions) also has growing
$\mathbb{S}_+^2$-rank.

\subsection{Limitations of $(\mathbb{S}_+^2)^p$-lifts}
\label{sec:js-fawzi}
We now state the main results of~\cite{fawzi2016representing}, establishing
fundamental limitations on the convex sets that can be exactly expressed using
LMIs with only $2\times 2$ blocks.
\begin{theorem}
	\label{thm:fawzi}
	There is no finite $p$ such that $PSD_{1,4}$, the cone of non-negative
	univariate quartic polynomials, has an $(\mathbb{S}_+^2)^p$-lift.
\end{theorem}
Since $PSD_{1,4} = SOS_{1,4}$ is the image of $\mathbb{S}_+^3$ under a linear
map (this follows from (\ref{eq:charac.sos.sdp}) as $Q$ is of size $3 \times 3$), the following is a direct consequence of
Theorem~\ref{thm:fawzi}.
\begin{corollary}[{\cite[Theorem 1]{fawzi2016representing}}]
	There is no finite $p$ such that $\mathbb{S}_+^3$, the cone of
	$3\times 3$ positive semidefinite matrices, has an $(\mathbb{S}_+^2)^p$-lift.
\end{corollary}
Providing the proof of these results is well beyond the scope of this article.
We will, however, sketch some of the ingredients.

First, we note that $PSD_{1,4}^* = \textup{cone}\{(1,t,t^2,t^3,t^4): t\in
\mathbb{R}\}$. This essentially follows directly from the definition of a
non-negative polynomial and the definition of the dual cone. Define a sequence
of points $\ell_j := (1,j,j^2,j^3,j^4)\in PSD_{1,4}^*$ for $j=1,2,\ldots$ and a
collection of non-negative quartic polynomials
\[ v_{\{i_1,i_2\}}(t) = \left[(i_1 - i_2)(i_1-t)(i_2-t)\right]^2,\]
indexed by pairs of positive integers $\{i_1,i_2\}$.  Then define a sequence of
$S^{(1)},S^{(2)},\ldots,S^{(k)},\ldots$ of $\binom{k}{2}\times k$ non-negative
matrices of the form
\[ S^{(k)}_{\{i_1,i_2\},j} = v_{\{i_1,i_2\}}(j) = \left[(i_1-i_2)(i_1-j)(i_2-j)\right]^2\]
where $1\leq i_1<i_2 \leq k$ and $1\leq j \leq k$. The aim is to show that the
$\mathbb{S}_+^2$-rank of the $S^{(k)}$ grows with $k$.  It then follows from
the discussion following the statement of Theorem~\ref{thm:psdrank} that
$PSD_{1,4} = SOS_{1,4}$ does not have a $(\mathbb{S}_+^2)^p$-lift for any
positive integer $p$.

Fawzi's argument makes crucial use of the sparsity pattern of these matrices,
and in particular certain relationships between the sparsity patterns
of $S^{(k)}$ and $S^{(k')}$ for different values of $k$ and $k'$.  In
particular, he shows that a certain combinatorial lower bound on the
$\mathbb{S}_+^2$-rank of $S^{(k)}$ must grow with $k$, and so the
$\mathbb{S}_+^2$-rank itself must grow with $k$.

\subsection{Challenges}
\label{sec:js-problems}
We conclude this section with a discussion of some challenges related to
understanding the limitations of approximating SOS cones with
convex cones having $(\mathbb{S}_+^2)^p$-lifts.

\subsubsection{Lower bounds on approximation quality}

Fawzi's result suggests that if we want to model SOS cones, in general, using
cones that have $(\mathbb{S}_+^2)^p$-lifts, approximation is necessary. One
approximation strategy uses SDSOS cones, but it is conceivable that a
significantly better approach exists. For a given notion of approximation and
approximation level $\epsilon$, it would be very interesting to produce lower
bounds on $p$, such that a given SOS cone can be $\epsilon$-approximated by a
convex cone having a $(\mathbb{S}_+^2)^p$-lift. Is there any reasonable notion
of approximation under which SDSOS is an optimal approximation in this sense?

\subsubsection{Strategies to construct approximate lifts}

On the positive side, how should we go about systematically approximating SOS
cones with convex cones having $(\mathbb{S}_+^2)^p$-lifts, with $p$ growing
mildly with approximation quality?  Are there ideas from classical constructive
approximation theory that can be  applied in this setting?  A concrete question
in this direction is the following:
\begin{quotation}
	For a given approximation quality $\epsilon$, how large must $p$ be so
	that we can $\epsilon$-approximate the cone of univariate non-negative
	polynomials of degree $d \geq 4$ with a convex cone having a
	$(\mathbb{S}_+^2)^p$-lift?
\end{quotation}
As an example of positive results in this broad direction, recent
work~\cite{fawzi2017semidefinite} develops an approach for constructing
high-quality approximations, with small $(\mathbb{S}_+^2)^p$-lifts, of
relative entropy cones.

\subsubsection{Techniques to lower bound $\mathbb{S}_+^2$-rank}
For a non-negative matrix $S$ there are numerous notions of cone rank. One
example is $\mathbb{S}_+^2$-rank, defined in Section~\ref{sec:js-bg} above;
another is the \emph{non-negative rank} (equivalent to $\mathbb{S}_+^1$-rank
in our notation); another is \emph{PSD rank} (see,
e.g.,~\cite{fawzi2015positive}, for a definition and survey related to its
properties).  These have interpretations in terms of modeling convex bodies in
terms of second-order cone programs, linear programs, and (general)
semidefinite programs, respectively. Recently, there has been considerable
interest, in a number of fields, in finding lower bounds on various notions of
cone rank such as these. Techniques for bounding the non-negative rank are most
developed, with lower bounds based on combinatorial
tools~\cite{fiorini2013combinatorial}, ideas from information
theory~\cite{braun2014information}, and systematic computational
methods~\cite{fawzi2016self}, being available. On the other end of the
spectrum, lower bounds on the PSD rank of non-negative matrices seem much more
challenging (see the survey~\cite{fawzi2015positive} for a discussion),
although there has been recent progress for some very specific matrices related
to combinatorial optimization problem~\cite{lee2015lower}.

The $\mathbb{S}_+^2$-rank in some ways behaves like the non-negative rank
(because of the inherent product structure), and in some ways like the PSD
rank. It would be very interesting to see which of the combinatorial tools
applicable to non-negative rank can be modified to this setting.  On the other
hand, developing methods for bounding the $\mathbb{S}_+^2$-rank that are
distinct from the approaches used for non-negative rank, may provide a path
towards understanding the PSD rank in general.

\balance
{\color{black}
\section{Conclusion}

In this paper, we have reviewed two new classes of techniques that aim to improve the scalability of SOS programs. First, two efficient first-order methods based on ADMM were proposed to solve the SDPs arising from SOS programs efficiently. Both of these techniques exploit the underlying sparsity to increase the computational efficiency.  
Second, we introduce techniques that replace the semidefinite program underlying any SOS program by more tractable convex programs such as linear or second-order cone programs. For this strategy, we first inner approximate the set of positive semidefinite matrices by the set of diagonally dominant matrices (resp. scaled diagonally dominant matrices) which is LP-representable (resp. SOCP representable). We then iteratively improve on these initial approximations while staying in the realm of LP and SOCP-representable sets. Finally, we reviewed recent results relating to how well one can represent the cone of SOS polynomials by SDPs involving only small blocks. We focused on second order representable cones, i.e., cones that involve SDPs with blocks of size bounded by 2. We presented a recent result that states that one cannot even represent the set of nonnegative univariate quartics using these kinds of cones. This leaves open the question of how well one can approximate SOS polynomials with SDPs involving only small blocks.
}


\bibliographystyle{IEEEtran}
\bibliography{BibCDC17}

\end{document}